\documentclass[12pt]{article}
\usepackage[final]{epsfig}
\usepackage{graphics}
\usepackage{amsmath}
\usepackage{amsfonts}
\usepackage{latexsym}
\usepackage{amssymb}
\usepackage{graphicx}
\usepackage{url}
\usepackage{epstopdf}

\newtheorem{lemma}{Lemma}[section]
\newtheorem{proposition}[lemma]{Proposition}
\newtheorem{remark}[lemma]{Remark}

\newtheorem{theorem}{Theorem}

\newtheorem{corollary}[lemma]{Corollary}

\newcommand{\proofend}{$\Box$\bigskip}

\newcommand{\R}{{\mathbb R}}

\def\proof{\paragraph{Proof.}}
\newcommand{\V}{{\rm Vol\ }}

\begin{document}

\title{Circumcenter of Mass and generalized Euler line}

\author{Serge Tabachnikov\footnote{
Department of Mathematics,
Pennsylvania State University, 
University Park, PA 16802, USA, 
tabachni@math.psu.edu} \and Emmanuel Tsukerman\footnote{Stanford University, emantsuk@stanford.edu}
}

\date{}
\maketitle

\begin{abstract}
We define and study a variant of the center of mass of a polygon and, more generally, of a simplicial polytope which we call the {\it Circumcenter of Mass (CCM)}.  The Circumcenter of Mass is an affine combination of the circumcenters of the simplices in a triangulation of a polytope, weighted by their volumes.  For an inscribed polytope, CCM coincides with the circumcenter. 

Our motivation comes from the study of completely integrable discrete dynamical systems, where the CCM is an invariant of the {\it discrete bicycle (Darboux) transformation} and of {\it recuttings} of polygons. 

We show that the CCM satisfies an analog of Archimedes' Lemma, a familiar property of the center of mass. We define and study a {\it generalized Euler line} associated to any simplicial polytope, extending the previously studied Euler line associated to the quadrilateral. We show that the generalized Euler line for polygons consists of all centers satisfying natural continuity and homogeneity assumptions and Archimedes' Lemma. 

Finally, we show that CCM can also be defined in the spherical and hyperbolic  settings.
\end{abstract}

\section{Introduction} \label{intro}

This note concerns  the following geometric construction. Let $P=(V_1,V_2,\dots,\\V_n)$ be an oriented $n$-gon in the plane and let $O$ be a point not on the sides of $P$ (or their extensions). One has a triangulation of $P$ by $n$ oriented triangles $OV_iV_{i+1}$ where the index is understood cyclically mod $n$; let $C_i$ be the circumcenter of the $i$th triangle and $A_i$ its signed area. Then the sum $A_1+\dots+A_n$ is the signed area $A(P)$ of the polygon; we assume that $A(P)\neq 0$. Consider the weighted sum of the circumcenters
\begin{equation} \label{defCCM}
CCM(P)=  \sum_{i=1}^n \frac{A_i}{A(P)}\ C_i.
\end{equation}
We call this point the {\it Circumcenter of Mass} of the polygon; we shall show that it does not depend on the choice of the point $O$.

\begin{figure}[hbtp]
\centering
\includegraphics[width=2.6in]{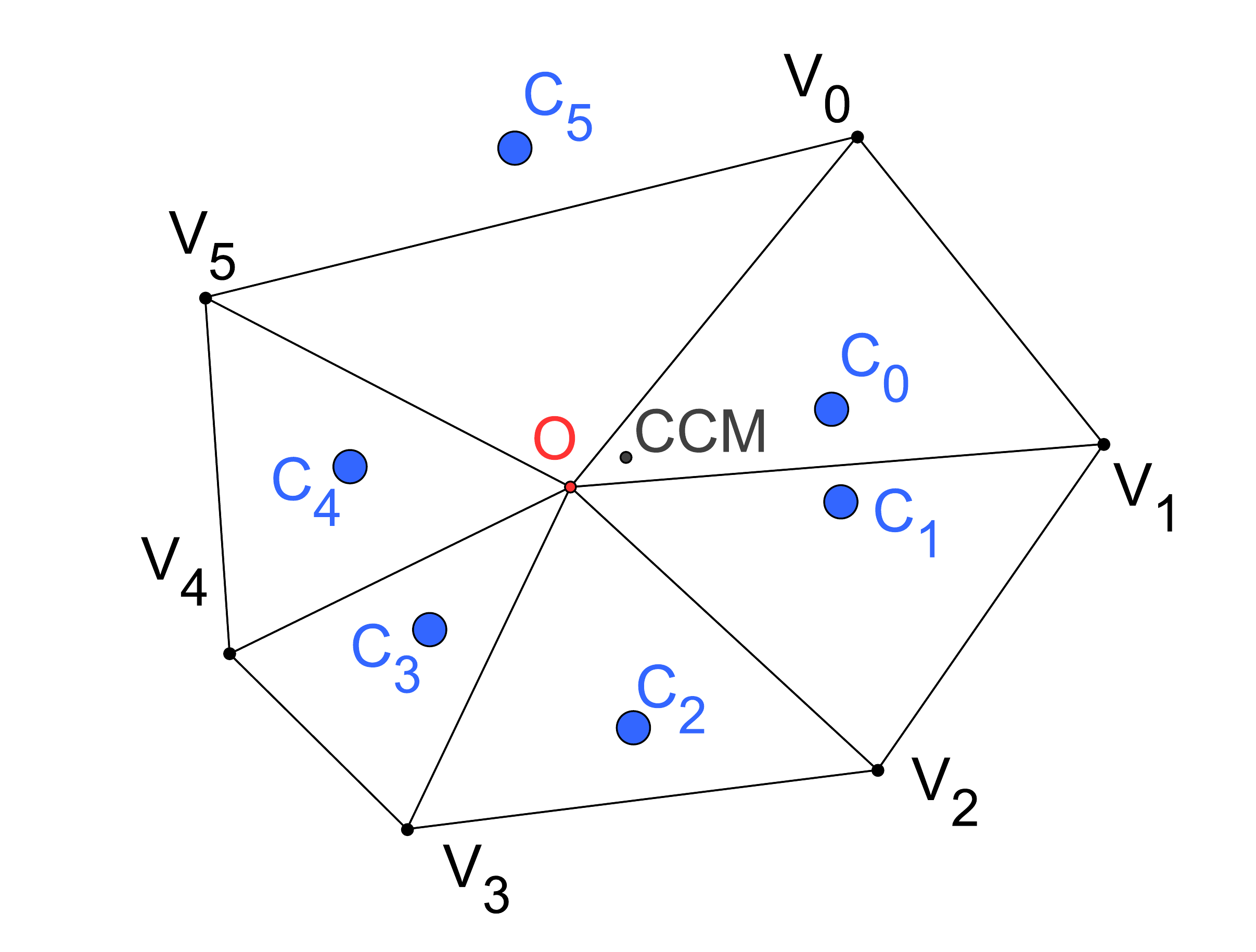}
\caption{Construction of the Circumcenter of Mass}
\label{Defn}
\end{figure}

In other words, we break $P$ into triangles $OV_{i}V_{i+1}$
and replace each such triangle with a point-mass equal to its area
and located at its circumcenter. Then $CCM(P)$ is the center of mass of
these point-masses, see Figure \ref{Defn}.

We learned about this construction from \cite{Ad1} where this center is mentioned as a point invariant under the  transformation of the polygon called {\it recutting}. In \cite{TT}, we proved that the Circumcenter of Mass is invariant under another map, called the {\it discrete bicycle} (or discrete Darboux) transformation; see \cite{Ho,PSW}. Although our motivation comes from the study of these discrete dynamical systems, we believe that the Circumcenter of Mass is an interesting geometric construction and it deserves attention on its own right.

If one replaces the circumcenter of a triangle by its centroid then, instead of $CCM(P)$, one  obtains the center of mass of a polygon $P$ (more precisely, of the homogeneous lamina bounded by $P$), which we denote by $CM(P)$.  

In the case of $CM(P)$, one may safely allow degenerate triangles, whose centroids are well defined and whose areas vanish, and which consequently do not contribute to the sum. However, not all degenerate triangles are allowed in the definition of the Circumcenter of Mass: although their areas vanish, their circumcenters may lie at infinity, in which case formula (\ref{defCCM}) fails to make sense. This discontinuity with respect to the location of point $O$ is a subtle point of the definition of $CCM(P)$.  

The contents of this paper are as follows. In Section \ref{welldef}, we prove that the definition of the Circumcenter of Mass is correct, that is, does not depend on the choice of the point $O$.

In Section \ref{prop}, we establish a number of properties of the Circumcenter of Mass:
\begin{enumerate}
\item The circumcenter of mass satisfies the Archimedes Lemma: if a polygon is a union of two
smaller polygons, then the circumcenter  of mass of the compound polygon is a weighted sum of the circumcenters of mass
of the two smaller polygons. As a consequence, one may use any triangulation of $P$ to define $CCM(P)$.
\item If $P$ is an equilateral polygon then $CCM(P)=CM(P)$.
\item One can naturally extend the notion of the Circumcenter of Mass to smooth curves. However, this continuous limit of the Circumcenter of Mass coincides with the center of mass of the lamina bounded by the curve.
\end{enumerate}

We start Section \ref{other} with the observation that one can replace the circumcenter or the centroid of a triangle by their affine combination and define a new ``center" of a polygon. The affine combinations of the circumcenter and the centroid of a triangle form the Euler line  of the triangle. This line contains a number of other interesting points, for example, the orthocenter. Applying our triangulation construction to a polygon $P$ yields  a line that we call the generalized Euler line  of $P$. 

We  investigate whether there are other natural assignments of ``centers" to polygons, satisfying  natural properties: the center depends analytically on the polygon, commutes with dilations, and satisfies the Archimedes Lemma. We show that all such centers are points of the generalized Euler line.

We also show that the Euler line is sensitive to symmetries, both Euclidean (reflections, rotations), and others, such as equilateralness.

In Section \ref{highdim}, we show that the construction of the Circumcenter of Mass extends to simplicial polyhedra in $\R^n$ in a rather straightforward manner. In Section \ref{sphere}, we extend our constructions to the spherical and hyperbolic geometries.

Whenever possible, we try to provide both geometric and algebraic proofs, so some statements in this paper are given multiple proofs. 
\medskip

{\it A historical note}. After this paper was finished, we learned from B. Gr\"unbaum that the notion of the circumcenter of mass appeared in the book by Laisant \cite{La} under the name ``pseudo-centre"; the discovery is credited to the Italian algebraic geometer Giusto Bellavitis. B. Gr\"unbaum and G. C. Shephard obtained some of the result that appear below in 1993, but they did not publish this material. We are grateful to B. Gr\"unbaum for this information.

\section{Correctness of definition. Degenerate triangles} \label{welldef}

Let us prove that the Circumcenter of Mass is well defined. 

We shall use the following notation for the CCM of
two polygons. If $P$ and $Q$ are polygons,  define 
$$
CCM(P\oplus Q)=\frac{A(P)\ CCM(P)+A(Q)\ CCM(Q)}{A(P)+A(Q)}.
$$
It is easy to check that this operation is commutative and associative. We also define
$$
CCM(P\ominus Q)=\frac{A(P)\ CCM(P)-A(Q)\ CCM(Q)}{A(P)-A(Q)}.
$$
We consider $CCM(\ominus P)$ to have the same CCM as $P$ but with mass of opposite sign.

\begin{lemma} \label{quad}
Given a quadrilateral  $ABCD$, one has:
$$
CCM(ABC\oplus CDA)=CCM(BCD\oplus DAB).
$$
\end{lemma}

\begin{figure}[hbtp]
\centering
\includegraphics[width=2.6in]{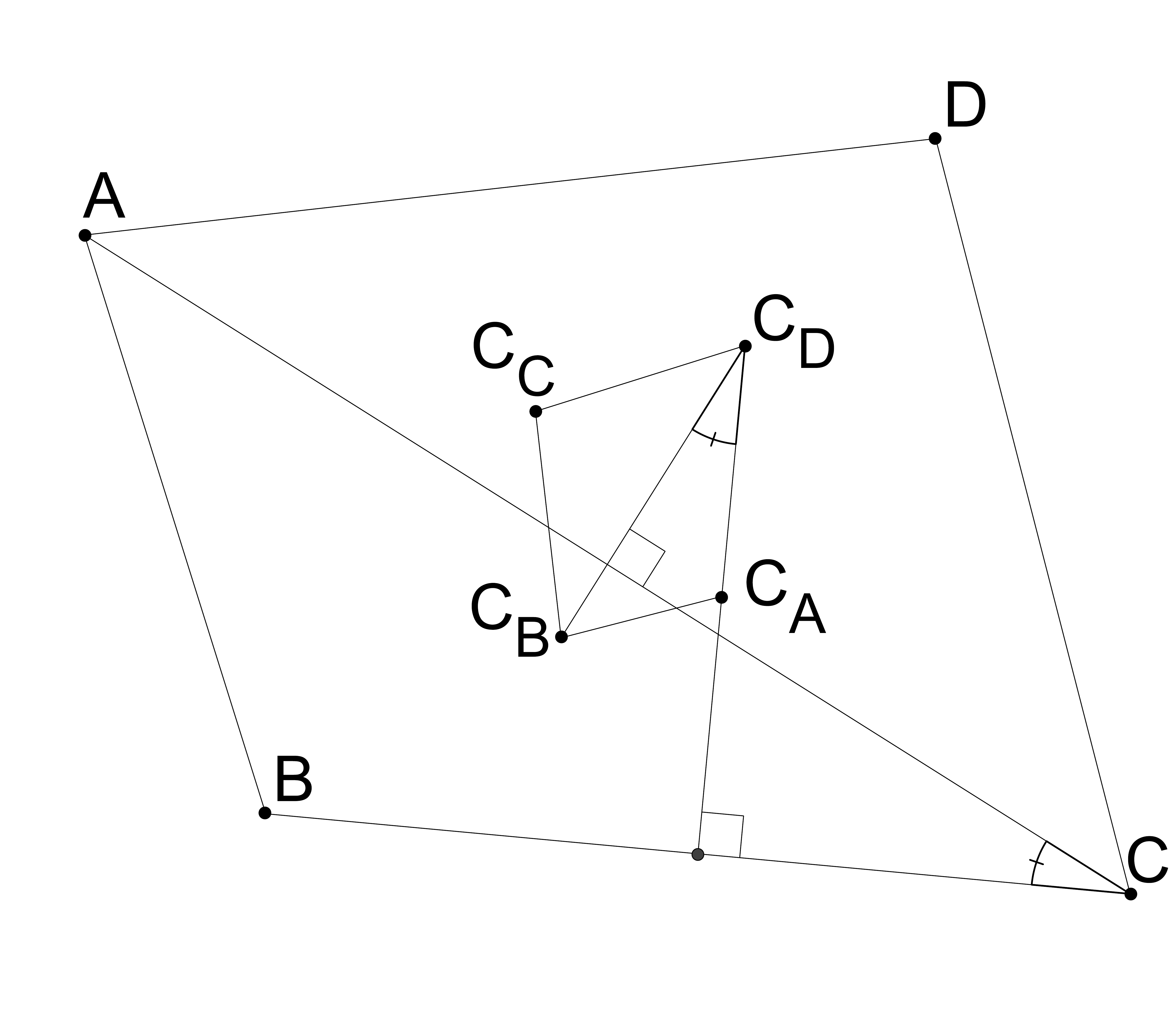}
\caption{Proof of Lemma \ref{quad}}
\label{orthologic}
\end{figure}

\proof 
Let $C_{A}$, $C_{B}$, $C_{C}$ and $C_{D}$ be the circumcenters
of $\triangle BCD$, $\triangle ACD$, etc., see Figure \ref{orthologic}. It suffices to see that $C_{A}C_{C}\cap C_{B}C_{D}$, divides
the diagonals of $C_{A}C_{B}C_{C}C_{D}$ at the ratio corresponding
to the areas of the triangles. We show this for one pair, and by symmetry, it will follow for the other.

An easy angle count shows that $\angle C_{B} C_{D} C_{D}=\angle BAC$, $\angle C_{A} C_{D} C_{B}=\angle A C B$, $\angle C_{C} C_{B} C_{D} = \angle C A D$, $\angle C_{A} C_{C} C_{D} =\angle A B D$, $\angle C_{C} C_{A} C_{D} = \angle C B D$. Indeed, the angles $C_{A} C_{D} C_{B}$ and $A C B$ have pairwise perpendicular sides, see Figure \ref{orthologic}, and likewise with the other pairs of angles.

To see that the $CCM$ divides diagonal $C_{A}C_{C}$ in ratio $A(ABD):A(CBD)$, it is enough to show that $A(ABD):A(CBD)=A(C_{D} C_{C} C_{B}):A(C_{D} C_{A} C_{B})$. Using the angles from above and the Law of Sines, 
\begin{equation*}
\begin{split}
\frac{A(C_{D} C_{C} C_{B})}{A(C_{D} C_{A} C_{B})}&=\frac{C_{C} C_{D} \cdot C_{B} C_{D} \sin(\angle BAC)}{C_{A} C_{D} \cdot C_{B} C_{D} \sin(\angle ACB)}=\frac{C_{C} C_{D} \sin(\angle BAC)}{C_{A} C_{D} \sin(\angle ACB)} 
=\frac{C_{C} C_{D} \cdot BC}{C_{A} C_{D} \cdot AB}\\
&=\frac{BC \sin(\angle CBD)}{AB \sin(\angle ABD)}=\frac{BC \cdot BD \sin(\angle CBD)}{AB \cdot BD \sin(\angle ABD)}=\frac{A(ABD)}{A(CBD)},
\end{split}
\end{equation*}
as needed.
\proofend

We use Lemma \ref{quad} to prove that the Circumcenter of Mass is well defined.  We wish to show that $CCM(\oplus_{i=0}^{n-1}V_{i}V_{i+1}O)$ is independent of $O$. We proceed by induction on $n$, the number
of vertices of the polygon $P$. 

The base case is $n=3$. By Lemma
\ref{quad}, 
\begin{equation*}
\begin{split}
&CCM(V_{0}V_{1}O\oplus V_{1}V_{2}O\oplus V_{2}V_{0}O)=CCM(V_{0}V_{1}O\ominus V_{2}V_{1}O\ominus OV_{0}V_{2})\\
&=CCM(V_{0}V_{1}O\ominus V_{1}OV_{0}\ominus V_{0}V_{2}V_{1})=CCM(\ominus V_{0}V_{2}V_{1})=CCM(V_{0}V_{1}V_{2}).
\end{split}
\end{equation*}

Now assume the inductive hypothesis for $n-1\geq3$ and consider the case for $n$:
\begin{equation*}
\begin{split}
CCM(\oplus_{i=0}^{n}V_{i}V_{i+1}O)&=CCM((\oplus_{i=0}^{n-2}V_{i}V_{i+1}O)\oplus V_{n-2}V_{n-1}O\oplus V_{n-1}V_{0}O)\\
&=CCM((\oplus_{i=0}^{n-2}V_{i}V_{i+1}O)\ominus V_{n-1}V_{n-2}O\ominus OV_{0}V_{n-1}).
\end{split}
\end{equation*}

By Lemma \ref{quad}, this expression is equal to
\begin{equation*}
\begin{split}
&CCM((\oplus_{i=0}^{n-2}V_{i}V_{i+1}O)\ominus V_{n-2}OV_{0}\ominus V_{0}V_{n-1}V_{n-2})\\
&=CCM((\oplus_{i=0}^{n-2}V_{i}V_{i+1}O)\oplus V_{n-2}V_{0}O)\oplus CCM(V_{n-2}V_{n-1}V_{0}).
\end{split}
\end{equation*}
The first term of the last expression is the CCM of a closed $(n-1)$-gon. By induction,
it is independent of $O$. The second term is clearly independent
of $O$. This completes the induction.
\medskip

Next, we deduce coordinate expressions for the Circumcenter of Mass.  Let $O$ be the origin and let $V_i=(x_i,y_i)$.

\begin{proposition} \label{coord}
One has: $CCM(P)=$
\begin{equation*}
\begin{split}
&\frac{1}{4 A(P)} \left(\sum_{i=0}^{n-1}-y_{i}y_{i+1}^{2}+y_{i}^{2}y_{i+1}+x_{i}^{2}y_{i+1}-x_{i+1}^{2}y_{i},\sum_{i=0}^{n-1}-x_{i+1}y_{i}^{2}+x_{i}y_{i+1}^{2}+x_{i}x_{i+1}^{2}-x_{i}^{2}x_{i+1}\right)\\
&=\frac{1}{4 A(P)}\left(\sum_{i=0}^{n-1}y_{i}(x_{i-1}^{2}+y_{i-1}^{2}-x_{i+1}^{2}-y_{i+1}^{2}),\sum_{i=0}^{n-1}-x_{i}(x_{i-1}^{2}+y_{i-1}^{2}-x_{i+1}^{2}-y_{i+1}^{2})\right).
\end{split}
\end{equation*}
\end{proposition}

\proof
The perpendicular bisector to $OV_{i}$ has equation
$$
y=-\frac{x_{i}}{y_{i}}\left(x-\frac{x_{i}}{2}\right)+\frac{y_{i}}{2}.
$$
Intersecting with the perpendicular bisector to $OV_{i+1}$, we find
that the coordinates of the circumcenter of $\triangle OV_{i}V_{i+1}$
are
$$
x=-\frac{x_{i+1}^{2}y_{i}-x_{i}^{2}y_{i+1}-y_{i}^{2}y_{i+1}+y_{i}y_{i+1}^{2}}{2(x_{i}y_{i+1}-x_{i+1}y_{i})}, \ 
y=-\frac{x_{i}^{2}x_{i+1}-x_{i}x_{i+1}^{2}+x_{i+1}y_{i}^{2}-x_{i}y_{i+1}^{2}}{2(x_{i}y_{i+1}-x_{i+1}y_{i})}.
$$
The denominators are four times the areas:
$$
2 A(\triangle OV_{i}V_{i+1}) = 
\det(\vec{OV_{i}},\vec{OV_{i+1}})=\det\left(\begin{array}{cc}
x_{i} & x_{i+1}\\
y_{i} & y_{i+1}
\end{array}\right)=x_{i}y_{i+1}-x_{i+1}y_{i}.
$$
This implies the result.
\proofend

\begin{remark} \label{turn}
{\rm  
Rotating the Circumcenter of Mass by $90^{\circ}$, we get
\begin{equation} \label{rot}
\frac{1}{4 A(P)}\sum_{i=0}^{n-1}|V_{i}|^{2}(V_{i+1}-V_{i-1}).
\end{equation}
This formula is given in \cite{Ad2} as an integral of the recutting of polygons. 

It is also interesting to compare with the formula for the center of mass of a homogeneous polygonal lamina; this formula has a similar structure (the ratio of a cubic polynomial and the area): $CM(P)=$
$$
\frac{1}{6 A(P)} \left(\sum_{i=0}^{n-1}(x_{i}+x_{i+1})(x_{i}y_{i+1}-x_{i+1}y_{i}),\sum_{i=0}^{n-1}(y_{i}+y_{i+1})(x_{i}y_{i+1}-x_{i+1}y_{i})\right).
$$
}
\end{remark}

We now use Proposition \ref{coord} to give a different proof of the independence of $CCM(P)$ of the choice of the origin $O$. 
We wish to show that parallel translation of the point $O$ through some vector results in translating $CCM(P)$ by the same vector.

Let the translation vector be $(\xi,\eta)$, that is,
$$
\bar x_i=x_i+\xi, \ \bar y_i = y_i + \eta.
$$
Let us compute the formula of Proposition \ref{coord} in terms of $(\bar x_i,\bar y_i)$. The area $A(P)$ is translation invariant.
 Let us consider the first coordinate, multiplied by $4A(P)$. As a polynomial in $\xi,\eta$, it is:
 $$
 A+B\xi+C\eta+D\xi\eta+E\eta^2,
 $$
 where $A,B,C,D,E$ are polynomials in $x_i,y_i$. The term $A$ is the same as in Proposition 2.3; as to the rest, one easily computes:
\begin{equation*}
\begin{split}
B=\sum 2(x_{i-1}-x_{i+1})y_i,\ &C=\sum [x_{i-1}^2-x_{i+1}^2+y_{i-1}^2-y_{i+1}^2 +2y_i(y_{i-1}-y_{i+1})],\\
&D= \sum 2(x_{i-1}-x_{i+1}),\ E= \sum 2(y_{i-1}-y_{i+1}),
\end{split}
\end{equation*}
where all sums are cyclic. We see that 
$$
B=4A(P),\ C=D=E=0,
$$
the latter equalities due to ``telescoping" of the cyclic sums. Therefore, the first coordinate of $CCM(P)$ is shifted by $\xi$. A similar computation shows that the second coordinate is shifted by $\eta$.
\medskip

\begin{figure}[hbtp]
\centering
\includegraphics[width=4in]{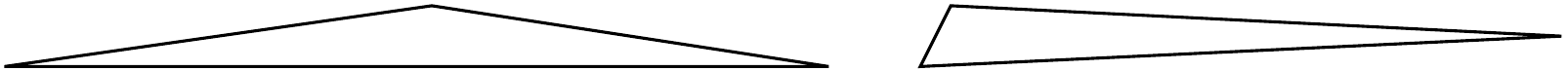}
\caption{The triangle on the left is dangerous and the one on the right is safe}
\label{degen}
\end{figure}

Next, we comment on the issue of degenerate triangles in the definition of the Circumcenter of Mass. A degenerate triangle is ``dangerous" if its circumcenter is at infinity, that is, if the circumradius is infinite. This radius is given by the formula
$$
R=\frac{a}{2\sin\alpha}=\frac{b}{2\sin\beta}=\frac{c}{2\sin\gamma},
$$
where $a,b,c$ are the side lengths of a triangle and $\alpha, \beta, \gamma$ are the respective angles. This makes it possible to distinguish between dangerous and safe degenerations, see Figure \ref{degen}: a degeneration is dangerous if an angle tends to zero or $\pi$ while the opposite side does not tend to zero. 

\begin{figure}[hbtp]
\centering
\includegraphics[width=2in]{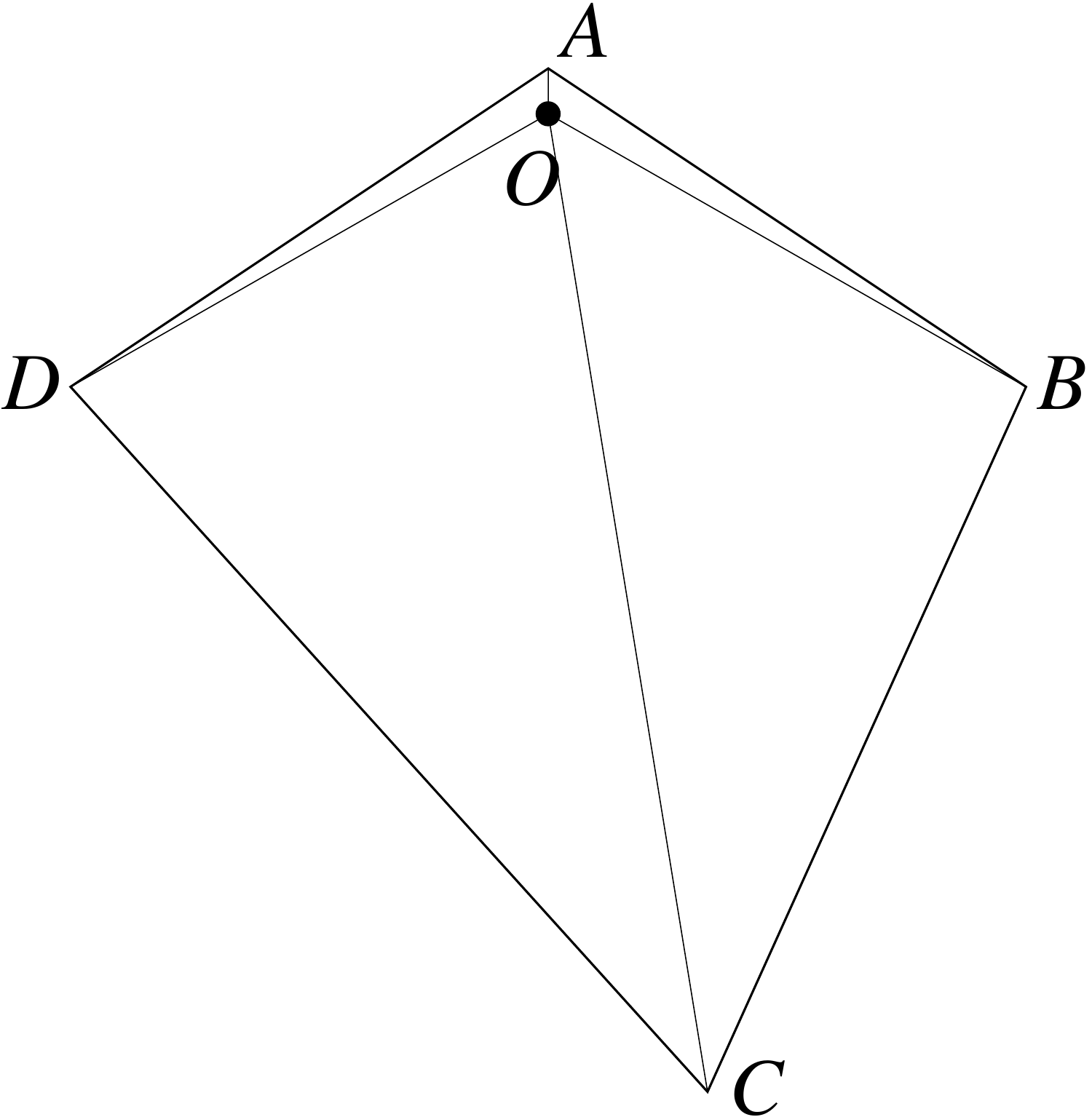}
\caption{An alternative proof of Lemma \ref{quad}}
\label{limit}
\end{figure}

\begin{remark} \label{altpf}
{\rm  Lemma \ref{quad} follows from the independence of $CCM(P)$ on the choice of the point $O$. Indeed, choosing the point $O$ close to vertex $A$ on the bisector of the angle $DAB$, see Figure \ref{limit}, one obtains, in the limit, two safe degenerate triangles $OAD$ and $OAB$. With this choice of $O$, one has $CCM(ABCD)=CCM(ABC\oplus CDA)$, and likewise for the the other pair of triangles.}
\end{remark}

\section{Properties of the Circumcenter of Mass.\\ Archimedes Lemma} \label{prop}

The center of mass satisfies the Archimedes Lemma which states the following: if an object is divided into two
smaller objects, the center of mass of the compound object lies on
the line segment joining the centers of mass of the two smaller objects; see, e.g., 
\cite{AM}. It turns out, the Circumcenter of Mass satisfies the Archimedes Lemma as well. 
 
\begin{figure}[hbtp]
\centering
\includegraphics[width=2in]{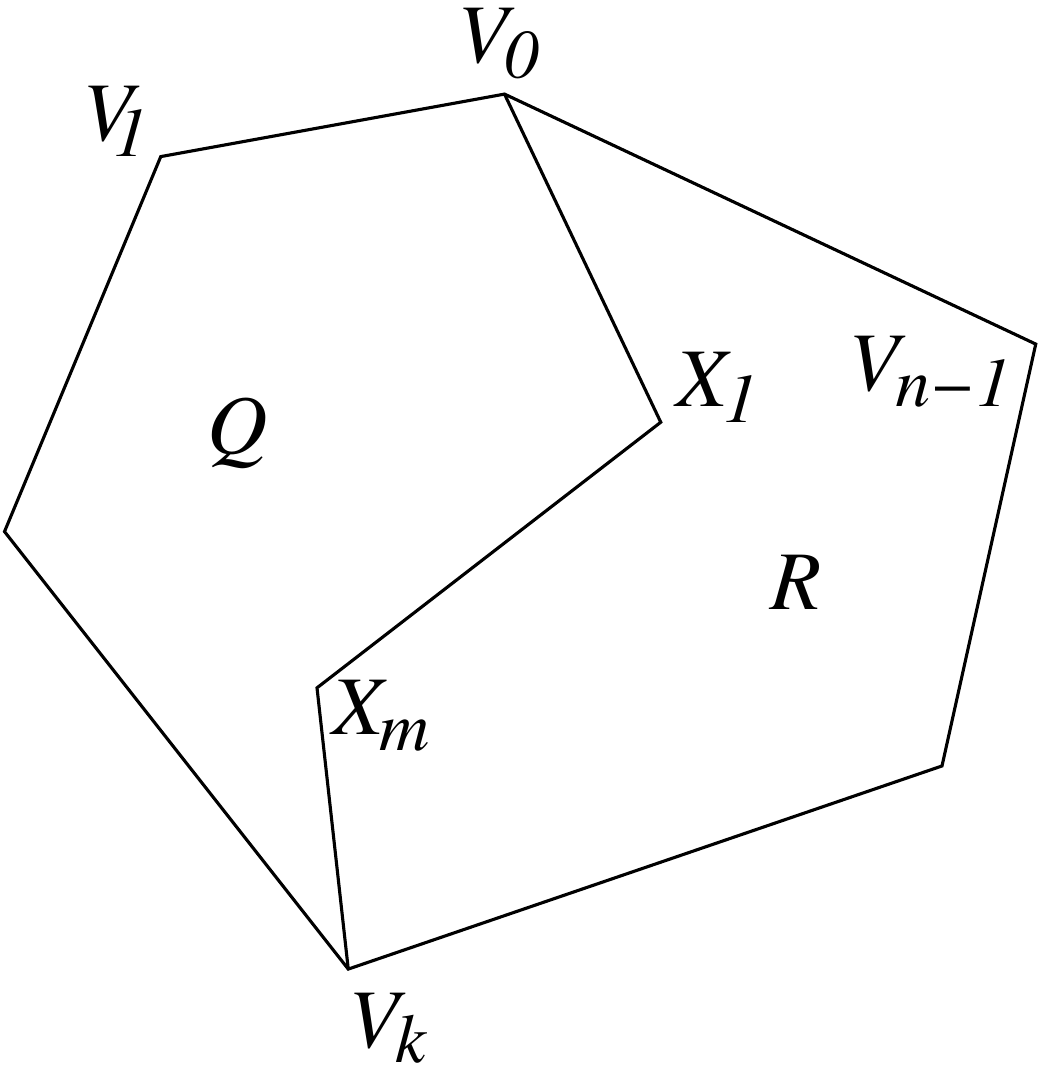}
\caption{Archimedes Lemma}
\label{cut}
\end{figure}

\begin{theorem} \label{Archimedes}
Let $P=V_{0}V_{1}\dots V_{n-1}$ be a  polygon, and let $V_{0}X_{1}X_{2}\dots X_{m}V_{k}$
be a polygonal line from vertex $V_{0}$ to  vertex $V_{k}$. 
Let $Q=V_{0}V_{1}\dots V_{k}X_{m}X_{m-1}\dots X_{1}$
and $R=V_{0}X_{1}X_{2}\dots X_{m}V_{k}V_{k+1}\dots V_{n-1}$ be the two closed polygons resulting from cutting $P$ along the polygonal
line. Then
$$
CCM(P)=CCM(Q\oplus R),
$$
see Figure \ref{cut}.
\end{theorem}

\proof
The idea of the proof is very simple: the contribution of the triangles whose sides constitute the cut cancels. More specifically, fix a generic point $O$. Then
$$
CCM(P)=\frac{1}{ A(P)}\sum_{i=0}^{n-1}A(\triangle OV_{i}V_{i+1})\ {CCM}(\triangle OV_{i}V_{i+1}).
$$
Also, $A(Q)\ CCM(Q)=$
\begin{equation*}
\begin{split}
 \left(\sum_{i=0}^{k-1}  A(\triangle OV_{i}V_{i+1})\  {CCM}(\triangle OV_{i}V_{i+1})\right)+
 A(\triangle OV_{k}X_{m})\ {CCM}(\triangle OV_{k}X_{m})\\
-\left(\sum_{i=1}^{m-1} A((\triangle OX_{i}X_{i+1}))\ {CCM}(\triangle OX_{i}X_{i+1})\right)+
A(\triangle OX_{1}V_{0})\ {CCM}(\triangle OX_{1}V_{0}),
\end{split}
\end{equation*} 
and $A(R)\ CCM(R)=$
\begin{equation*}
\begin{split}
 A(\triangle OV_{0}X_{1})\  {CCM}(\triangle OV_{0}X_{1})+\left(\sum_{i=1}^{m-1}A(\triangle OX_{i}X_{i+1})\ {CCM}(\triangle OX_{i}X_{i+1})\right)\\
+A(\triangle OX_{m}V_{k})\ {CCM}(\triangle OX_{m}V_{k})+\left(\sum_{i=k}^{n-1}A(\triangle OV_{i}V_{i+1})\  {CCM}(\triangle OV_{i}V_{i+1})\right).
\end{split}
\end{equation*}
Adding and canceling terms,  the desired equality follows.
\proofend

\begin{remark}
{\rm One subtelty of Archimedes Lemma has to do with degenerate triangulations. At first sight, a counterexample to Theorem \ref{Archimedes} is a right isosceles triangle $P$, which is bisected into two smaller right isosceles triangles $Q$ and $R$ along its axis of symmetry. In this case CCM(P) is the midpoint of the base of P, which is not on the line joining CCM(Q) and CCM(R), the midpoints of its other two sides. \\
\indent The issue is that there is a hidden degenerate triangle. Label the isosceles triangle ABC with B having a right angle and let D be the foot of the perpendicular from B to side AC. If we perturb D away from line AC, we obtain an additional triangle ACD. This triangle has minute area but its circumcenter is far away. In this situation everything works out as expected. The problem occurs when we disregard this triangle because of its zero area. In other words, the CCM of ABC is the weighted average of the CCMs of ABD, BCD and ACD.}
\end{remark}

As a consequence of Theorem \ref{Archimedes}, one can use any triangulation of $P$ to define $CCM(P)$, not only a triangulation obtained by connecting point $O$ to the vertices. In particular, we have the following, expected, corollary.

\begin{corollary} \label{inscr}
The Circumcenter of Mass of an inscribed polygon is the circumcenter.
\end{corollary}

\proof
Let $P$ be inscribed, and let $O$ be the circumcenter. Consider a triangulation by diagonals. Then $O$ is the circumcenter of each triangle involved, and hence the Circumcenter of Mass is $O$.

Alternatively, one can use formula (\ref{rot}). For a circumscribed polygon, all $|V_i|$ are equal, hence (\ref{rot}) yields zero, that is, the origin $O$. Rotating back $90^{\circ}$ is $O$ as well.
\proofend

\begin{figure}[hbtp]
\centering
\includegraphics[width=1.2in]{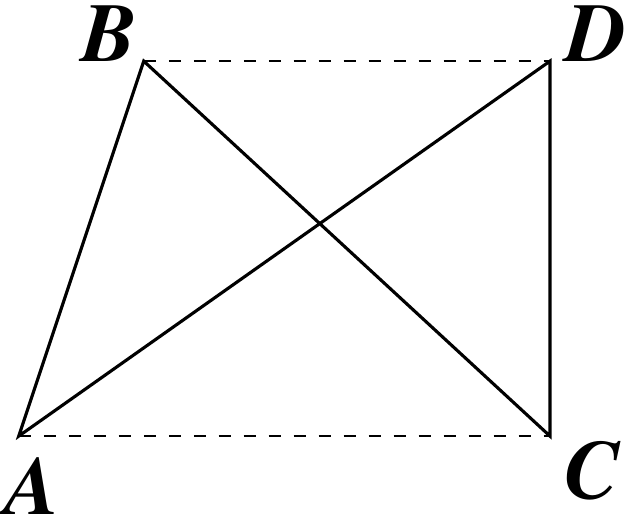}
\caption{The perpendicular bisectors of the diagonals of the quadrilateral $ABCD$ are parallel if and only if its signed area vanishes}
\label{zero}
\end{figure}

\begin{remark} \label{zeroarea}
{\rm It follows from Theorem \ref{Archimedes} that the Circumcenter of Mass of a quadrilateral is the intersection point of the perpendicular bisectors of its diagonals. These perpendicular bisectors may be parallel; this happens exactly when the area of the quadrilateral is equal to zero, see Figure \ref{zero}.
}
\end{remark}

Next, we consider the case when $P$ is equilateral, that is, all the sides have equal lengths.

\begin{theorem} \label{equilat}
If $P$ is an equilateral polygon then $CCM(P)=CM(P)$.
\end{theorem}

\proof
We present a computational proof based on Proposition \ref{coord} and Remark \ref{turn}. It would be interesting to find a geometric argument as well.

Using the formulas for $CCM(P)$ and $CM(P)$ from Proposition \ref{coord} and Remark \ref{turn}, we find:
\begin{equation} \label{2-3}
\begin{split}
2 CCM(P) - 3 CM(P) = \frac{1}{2A(P)}\ \Big(\sum (y_i^2y_{i+1}-y_iy_{i+1}^2-x_ix_{i+1}y_{i+1}+x_ix_{i+1}y_i),\\
\sum (x_ix_{i+1}^2-x_i^2x_{i+1}-x_iy_iy_{i+1}+x_{i+1}y_iy_{i+1} )\Big)\\
=\frac{1}{2A(P)}\ \Big(\sum (x_ix_{i+1}+y_iy_{i+1})(y_i-y_{i+1}), \sum (x_ix_{i+1}+y_iy_{i+1})(x_{i+1}-x_i)\Big). 
\end{split}
\end{equation}

Assume now that $P$ is equilateral:
$
(x_i-x_{i+1})^2+(y_i-y_{i+1})^2=1
$
for all $i$. Then
$$
x_i^2+y_i^2+x_{i+1}^2+y_{i+1}^2=1+2(x_ix_{i+1}+y_iy_{i+1}),
$$
and hence 
$$
x_{i-1}^2+y_{i-1}^2-x_{i+1}^2-y_{i+1}^2=2(x_{i-1}x_{i}+y_{i-1}y_{i}-x_{i}x_{i+1}-y_{i}y_{i+1}).
$$
Substitute this to the second formula of Proposition \ref{coord} to obtain
\begin{equation*}
\begin{split}
CCM(P)=\frac{1}{2A(P)}\ \big(\sum (y_ix_{i-1}x_i+y_{i-1}y_i^2-y_ix_ix_{i+1}-y_i^2y_{i+1}),\\
\sum (-x_{i-1}x_i^2-x_iy_{i-1}y_i+x_i^2x_{i+1}+x_iy_iy_{i+1})\big)\\
=\frac{1}{2A(P)}\ \big(\sum (x_ix_{i+1}+y_iy_{i+1})(y_{i+1}-y_i), \sum (x_ix_{i+1}+y_iy_{i+1})(x_i-x_{i+1})\big).
\end{split}
\end{equation*}
Comparing with (\ref{2-3}), we see that $CCM(P)= 3 CM(P)-2CCM(P)$, hence $CCM(P)=CM(P)$.
\proofend

Another case when the Circumcenter of Mass coincides with the center of mass is the continuous version of the former. Let $\gamma(t)$ be a parameterized smooth curve and $O$ be a point not on any tangent line to $\gamma$. In other words, $\gamma$ is star-shaped with respect to $O$. One defines the Circumcenter of Mass: 
$$
CCM(\gamma)=\frac{\int {C}(t)\ dA}{\int dA},
$$
where $C(t)$ denotes the limiting $\varepsilon \to 0$ position of the vector from $O$ to the circumcenter of the infinitesimal triangle $O \gamma(t) \gamma(t+\varepsilon)$, and $dA$ is the area of this infinitesimal triangle. Denote by $CM(\gamma)$ the center of mass of the homogeneous lamina bounded by $\gamma$.

\begin{theorem} \label{cont}
One has: $CCM(\gamma)=CM(\gamma)$.
\end{theorem}

\proof It is easiest to consider the continuous limit of formulas (\ref{rot}). Let $\gamma(t)=(x(t),y(t))$. Then the continuous limit of $V_{i+1}-V_{i-1}$ is $2\gamma'(t)$, and we obtain
$$
\frac{1}{2A(\gamma)} \left(\int (x^2+y^2)x' dt, \int (x^2+y^2)y' dt)  \right) = 
\frac{1}{A(\gamma)} \left(-\iint y dA, \iint x dA  \right), 
$$
where $dA=dx\wedge dy$ and the equality is due to the Stokes theorem. It remains to turn this vector back by $90^{\circ}$ to obtain
$$
\frac{1}{A(\gamma)} \left(\iint x dA, \iint y dA  \right), 
$$
 a well known formula for the center of mass.
\proofend

Thus one may view the Circumcenter of Mass of a polygon as a ``different" discretization of the center of mass of a lamina bounded by a continuous curve. 

\section{Generalized Euler line. Other centers} \label{other} 

The line through the center of mass and the circumcenter of a triangle is called the Euler line. The points of this line are affine combinations $C_t:= t CM + (1-t) CCM$; thus $C_0=CCM$ and $C_1=CM$. 
For example, the orthocenter of a triangle lies on the Euler line and is given by  $3CM-2CCM$, that is, the orthocenter is the point $C_{3}$.

For a fixed $t$, one can repeat the construction of the Circumcenter of Mass of a polygon to obtain the center $C_t(P)$: triangulate $P$ and take the sum of the centers $C_t$ of the triangles, weighted by their areas. The resulting center is again independent of the triangulation, and all these centers lie on a line that we call the {\it generalized Euler line} of the polygon $P$. 
Since the ratios of distances on the Euler line of a triangle are transferred over to those on the generalized Euler line, we conclude that $C_t(P) = t CM(P) + (1-t) CCM(P)$ for every polygon $P$.

We note that, for the case of quadrilaterals, the generalized Euler line was discussed in \cite{M}.
One remark we may add to the discussion is the next result, which follows easily from our approach.

\begin{figure}[hbtp]
\centering
\includegraphics[width=2.8in]{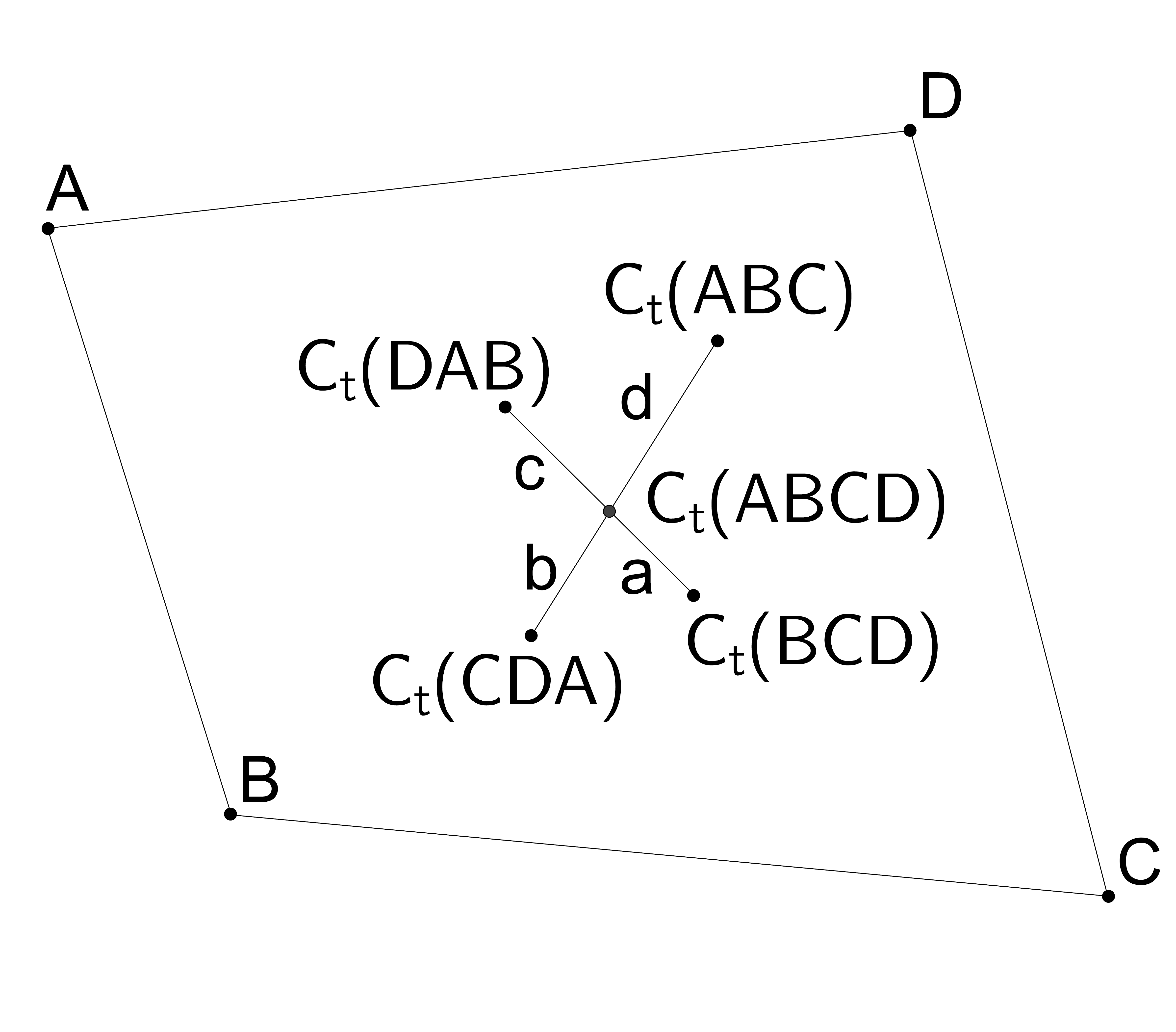}
\caption{Proposition \ref{ratios}}
\label{ratioTriangles}
\end{figure}

\begin{proposition} \label{ratios}
Let $a,b,c,d$ be the distances from $C_t(ABCD)$ to $C_t(BCD)$, $C_t(CDA)$, $C_t(DAB)$ and $C_t(ABC)$, respectively. Then 
$$
a A(BCD) = c A(DAB),  b A(CDA) = d A(ABC),
$$
see Figure \ref{ratioTriangles}.
\end{proposition}

Next, we shall show that the above described centers $C_t(P)$ are the only ones satisfying natural assumptions. Namely, assume that one assigns a ``center" to every polygon so that the center depends analytically on the polygon, commutes with dilations, and satisfies the Archimedes Lemma.

\begin{theorem} \label{unique}
A center satisfying the above assumptions is $C_t(P)$ for some $t$.
\end{theorem}

\proof 
Consider an isosceles triangle with base of length $2$ and base angles $\alpha$. By symmetry, the center is on the
axis of symmetry. Let it be at height $f(\alpha)$ above the base.
For example, if the center is $CM$ then $f(\alpha)={(\tan\alpha)}/{3}$,
and if the center is $CCM$ then $f(\alpha)=-\cot2\alpha$. 
From the symmetry of the equilateral triangle, $f({\pi}/{3})={1}/{\sqrt{3}}.$

Any triangle can be triangulated into three isosceles triangles, hence the function $f(\alpha)$ determines the center uniquely. 
We shall prove that $f$ satisfies a certain functional equation, solve it, and deduce the desired result.

\begin{figure}[hbtp]
\centering
\includegraphics[width=1in]{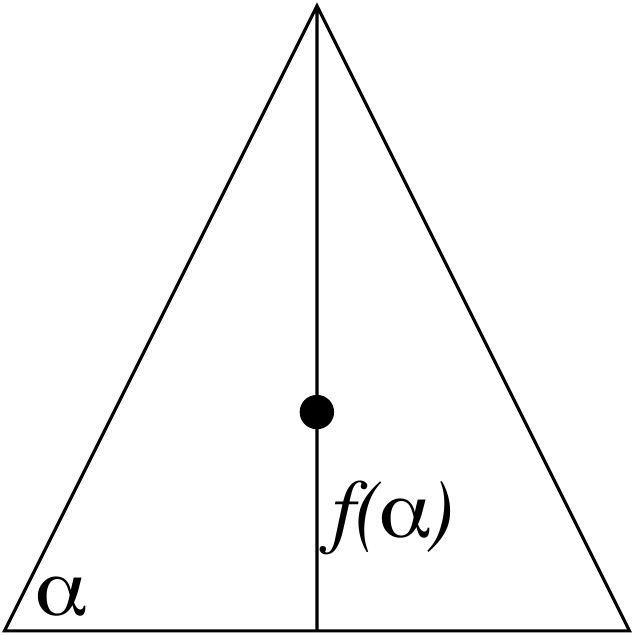}
\quad\quad\quad
\includegraphics[width=1in]{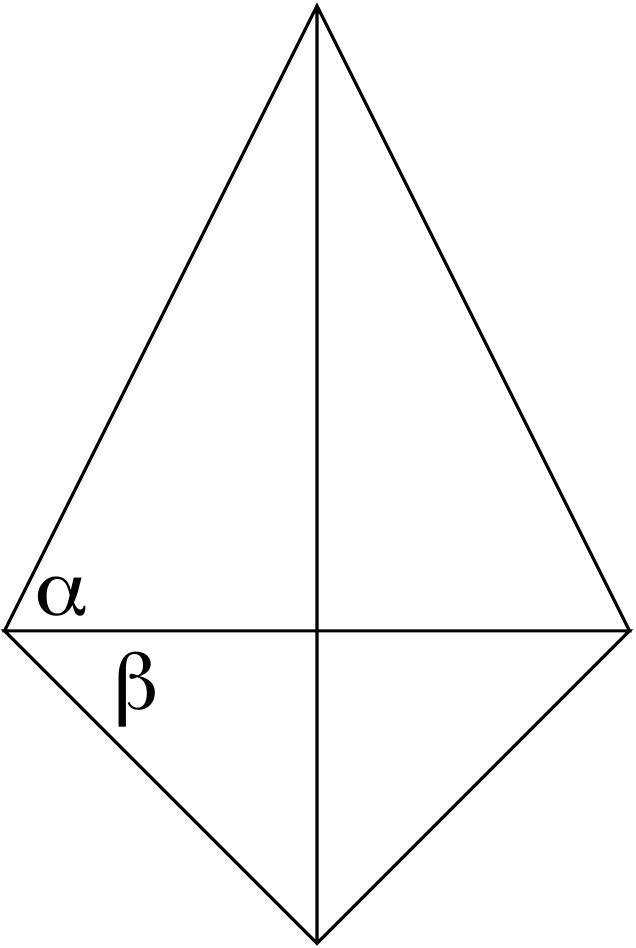}
\caption{Triangle and kite}
\label{triangle}
\end{figure}

Consider a kite made of two isosceles triangles with angles $\alpha$
and $\beta$. By the Archimedes Lemma, the center is  on the axis of symmetry at (signed) height
\begin{equation} \label{dist}
\frac{f(\alpha)\tan\alpha-f(\beta)\tan\beta}{\tan\alpha+\tan\beta},
\end{equation}
see Figure \ref{triangle}.

\begin{figure}[hbtp]
\centering
\includegraphics[width=2in]{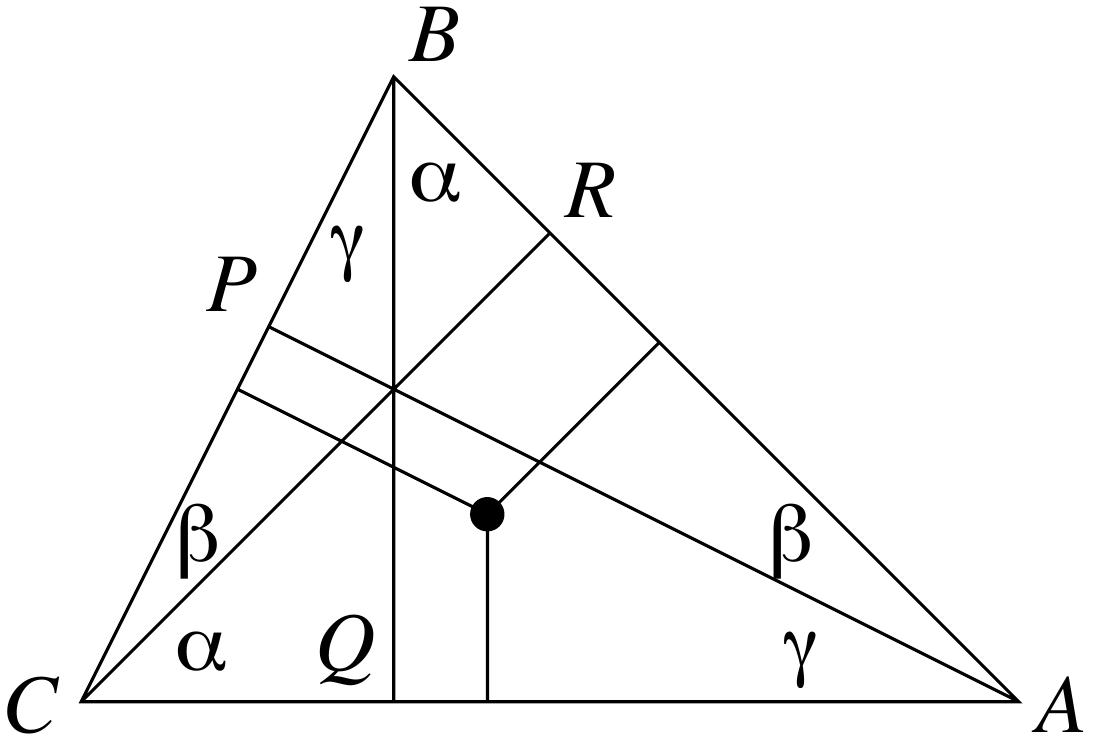}
\caption{A center of a triangle}
\label{reflections}
\end{figure}

Consider a triangle and its center, see Figure \ref{reflections}. One can reflect the triangle in either of its sides to obtain a kite. In this way, we can determine the signed distance from the center to each of the altitudes. For example, reflecting in the side $AC$, we use the (scaled) formula (\ref{dist}) to find that the distance to the altitude from vertex $B$ equals
$$
h_B \frac{f(\alpha)\tan\alpha-f(\gamma)\tan\gamma}{\tan\alpha+\tan\gamma},
$$
where $h_B=|BQ|$ is the length of this altitude.

Let $d_A, d_B, d_C$ be the signed distances from a point to the altitudes $AP, BQ,\\ CR$, respectively. 
The sign is determined by orienting the altitudes from the vertices of the triangle. 
We claim that
\begin{equation} \label{vanish}
\frac{d_A}{h_A} + \frac{d_B}{h_B} + \frac{d_C}{h_C} =0.
\end{equation}
Indeed, the left hand side of (\ref{vanish}) is a linear function on the plane, and it suffices to show that it vanishes in the vertices of the triangle. Consider vertex $A$. The triangles $ABQ$ and $ACR$ are similar, hence $|AQ|/|BQ|=|AR|/|CR|$. Taking the orientation into account, this proves (\ref{vanish}) for point $A$, and likewise for $B$ and $C$.

Thus (\ref{vanish}) implies the following functional equation on $f$:
\begin{equation} \label{main}
\frac{f(\alpha)\tan\alpha-f(\beta)\tan\beta}{\tan\alpha+\tan\beta}+\frac{f(\beta)\tan\beta-f(\gamma)\tan\gamma}{\tan\beta+\tan\gamma}+\frac{f(\gamma)\tan\gamma-f(\alpha)\tan\alpha}{\tan\gamma+\tan\alpha}=0
\end{equation}
for all triples of angles satisfying $\alpha+\beta+\gamma=\frac{\pi}{2}$. One can check that (\ref{main}) holds for $f(\alpha)={\tan\alpha}$ and for $f(\alpha)=\cot 2\alpha$. 

To solve equation (\ref{main}), we set $x=\tan\alpha$, $y=\tan\beta$, $z=\tan\gamma$
and $g(x)=xf(\arctan x).$ Then (\ref{main}) becomes
\begin{equation} \label{maing}
\frac{g(x)-g(y)}{x+y}+\frac{g(y)-g(z)}{y+z}+\frac{g(z)-g(x)}{z+x}=0
\end{equation}
subject to $\arctan x+\arctan y+\arctan z=\frac{\pi}{2}$. We have
\begin{equation} \label{elim}
z=\tan(\frac{\pi}{2}-\arctan x+\arctan y)=\frac{1-xy}{x+y}.
\end{equation}

We set $y=0$ in equation (\ref{maing}) to obtain
$$
\frac{g(x)-g(0)}{x}+\frac{g(0)-g(\frac{1}{x})}{\frac{1}{x}}+\frac{g(\frac{1}{x})-g(x)}{\frac{1}{x}+x}=0,
$$
or
$$
g(x)-x^4 g(\frac{1}{x})+g(0)(x^{4}-1)=0.
$$

Expand $g(x)$ in a power series: $g(x)=a_{0}+a_{1}x+a_{2}x^{2}+...$.
Then, equating terms of the same degrees, we find that $a_{n}=0$ for $n\geq5$, and $ a_4=0,  a_3=a_1$.

So $g(x)$ is of the form $a_{0}+a_{1}x+a_{2}x^{2}+a_{1}x^{3}$. Since
equation (\ref{maing}) is linear in $g$, it suffices
to consider each monomial individually. The $0$th, $1$st and $2$nd-degree
monomials satisfy the equation. For the $3$rd-degree, we have
$$
a_{1}\left(\frac{x^{3}-y^{3}}{x+y}+\frac{y^{3}-z^{3}}{y+z}+\frac{z^{3}-x^{3}}{x+z}\right)=0.
$$

If $a_{1}\neq0$, it is necessary for the expression in the brackets
to vanish identically after substituting for $z$ via equation (\ref{elim}).
However, it is easy to see that it does not by, e.g., substituting
$x=1$, $y=2$ and $z=-{1}/{3}$. It follows that $a_{1}=0$,
so that $g(x)=a_{0}+a_{2}x^{2}$. 

Therefore 
$
f(\alpha)=c_{1}\cot\alpha+c_{2}\tan\alpha.
$
Since $f({\pi}/{3})={1}/{\sqrt{3}}$, 
$$
\frac{c_{1}}{\sqrt{3}}+c_{2}\sqrt{3}=\frac{1}{\sqrt{3}} \implies c_{1}=1-3c_{2}.
$$
Hence 
$$
f(\alpha)=(1-3c)\cot\alpha+c\tan \alpha = t \left(\frac{\tan\alpha}{3}\right) + (1-t) (-\cot 2\alpha)
$$
with $t=3-6c$. Thus the solution set is the Euler line, as claimed.
\proofend

We now show that the generalized Euler line is highly sensitive to symmetries. As such, one might consider utilizing the Euler line to detect symmetries, or as a certificate for their non-existence. The symmetries involved are of different type, some being Euclidean and some of another sort, such as equilateralness.  

As we saw in the previous section from Theorem \ref{equilat}, if $P$ is an equilateral polygon, then the Euler line of $P$ degenerates to a point.

\begin{theorem} 
Let $P$ be a polygon  of nonzero area so that the Euler line $E$ is defined.
\begin{enumerate}
\item If $P$ has a line of reflection symmetry $L$, then $E$ is either $L$ or a point on $L$.
\item If $P$ has a center $C$ of rotational symmetry and none of the extensions of sides of $P$ pass through the center, then $E=C$.
\item Assume that the sides $S_1,S_2 \ldots, S_n$ of $P$ satisfy $|S_1|=|S_2|=\ldots=|S_{n-1}|$. Then $E$ is orthogonal to side $S_n$.
\end{enumerate}
\end{theorem}

\proof
(1). First, we note that $P$ cannot have sides which are subsets of $L$ because then at least one such side would have three edges incident to a vertex. Thus we may place $O$ on $L$ and triangulate. Since every triangle on one side of $L$ has a mirror image counterpart on the other side of $L$, the Circumcenter of Mass lies on $L$ and so does the center of mass. \\
(2) We let $O=C$ and triangulate from $O$. To each triangle in the triangulation, there is a class of  triangles obtained from it via rotation. The vector sum of the centers of mass of the triangles of this class is equal to $O$. \\
(3) Let $V_1,V_2,\ldots,V_n$ be the vertices of $P$. We reflect $P$ in side $S_n$. Denote the image of a vertex under reflection by a prime. Consider the polygon $Q=V_1 V_2 \cdots V_n V_{n-1}' \cdots V_{2}'$. This polygon has twice the area of $P$ and since it is equilateral, its Euler line is a point (the center of mass). Since $Q$ has $S_n$ as an axis of symmetry, its center of mass lies on $S_n$ (or its extension). By Archimedes' Lemma (Theorem \ref{Archimedes}), the Euler line of $Q$ is the weighted sum of the Euler lines of $P$ and its reflection. Since the area of the reflection is the same as that of $P$, the weights are the same. Let $e$ and $e'$ denote the two Euler lines. These are symmetric about $S_n$, so the Euler line of $Q$ must be the midpoint of every pair of mirror points of $e$ and $e'$. This is only possible if $e=e'$ and are orthogonal to $S_n$. 
\proofend

\begin{remark} \label{otherlines}
{\rm The generalized Euler line that we have defined  is as a property of the polygon, rather than the set of its vertices. For example, if a different choice of edges is given to a quadrilateral, e.g., rather than considering $ABCD$, we consider $ABDC$, its center of mass and Circumcenter of Mass will be different, and consequently its Euler line also.

In the case of a quadrilateral, a definition of the Euler line which is a property of the set of vertices is given in \cite{Ma}, following the work \cite{RT}. The authors of \cite{RT} construct a point which is a ``replacement" to the circumcenter of a quadrilateral when the quadrilateral is no longer cyclic. This point, called the ``isoptic point", is a center of the quadrilateral which depends only on the set of vertices, in the sense that it is independent of the choice of edges for the four vertices of the quadrilateral, e.g., it is the same for $ABCD$ and $ABDC$. The Euler line is then defined as the line connecting the isoptic point and the centroid of the vertices of the quadrilateral.  In \cite{Ts}, the second author extends the definition of the isoptic point to $(n+2)$-polytopes in $\mathbb{R}^n$, thus constructing an Euler line which depends only on the set of vertices of the $(n+2)$-tope.}
\end{remark}

\section{In higher dimensions} \label{highdim}

In this section, we extend the construction of the Circumcenter of Mass to simplicial polyhedra in $\R^n$. The construction is the same: given an oriented simplicial polyhedron $P$, choose a point $O$ that does not belong to any of the hyperplanes of its facet, and triangulate $P$ by the simplices with vertex $O$ and the bases the facets of $P$. Let $C_i$ be the circumcenters of these simplices, and define the circumcenter
\begin{equation*} \label{defhigh}
CCM(P)=  \sum \frac{V_i}{V(P)}\ C_i,
\end{equation*}
where $V_i$ denotes the signed volume of $i$th simplex, $V(P)$ is the volume of $P$, and the sum is taken over all facets of $P$.

Similarly to the 2-dimensional case, we have the following result.

\begin{theorem} \label{indep}
The Circumcenter of Mass is well defined, that is, does not depend on the choice of point $O$.
\end{theorem}

As a preparation to the proof, we state two lemmas. The first one is well known and we do not prove it.
Let $S=(V_0,V_1,\dots,V_n)$ be a simplex in $\R^n$. 

\begin{lemma} \label{vol}
The signed volume $\V (S)$ is 
$$
\frac{\det M(S)}{n!},
$$
 where $M$ is an $n+1$ by $n+1$ matrix
$$
\left(\begin{array}{cccc}
V_0&V_1&\dots&V_n\\
1&1&\dots&1
\end{array}\right) 
$$
\end{lemma}

The second lemma is lesser known so, for completeness, we provide a proof. Let $S=(V_0,V_1,\dots,V_n)$ be a simplex and $C(S)$ its circumcenter. 

\begin{lemma} \label{cc}
The $i$th coordinate of $C(S)$ is 
$$
\frac{\det M_i(S)}{2 \det M(S)},
$$
where the matrix $M_i$ is obtained from the matrix $M$ by replacing its $i$th row by the row $(|V_0|^2, |V_1|^2,\dots,|V_n|^2)$.
\end{lemma}

\proof Let $V_i=(x_{i1},\dots,x_{in}), \ i=0,\dots,n.$ Consider the equation
$$
\det \left|\begin{array}{cccccc}
x_1^2+\dots+x_n^2&x_{01}^2+\dots+x_{0n}^2&x_{11}^2+\dots+x_{1n}^2&\dots&x_{n1}^2+\dots+x_{nn}^2\\
x_1&x_{01}&x_{11}&\dots&x_{n1}\\
x_2&x_{02}&x_{12}&\dots&x_{n2}\\
\dots&\dots&\dots&\dots&\dots\\
x_n&x_{0n}&x_{1n}&\dots&x_{nn}\\
1&1&1&1&1
\end{array}\right|=0
$$
 in variables $x_1,\dots,x_n$. This  is an equation of a sphere, and this sphere passes through points $V_0,\dots,V_n$. Thus this is the equation of the circumsphere of the simplex $S$. 

Expand the determinant in the first column:
$$
(x_1^2+\dots+x_n^2)\det M  - x_1\det M_1 - x_2\det M_2 -\dots - x_n\det M_n +C=0,
$$
or 
$$
\left(x_1-\frac{\det M_1}{2\det M}\right)^2 + \dots + \left(x_n-\frac{\det M_n}{2\det M}\right)^2 = C_1,
$$
where $C,C_1$ are constants depending on the coordinates of the vertices of the simplex.
This implies the statement of the lemma.
\proofend

\paragraph{Proof of Theorem \ref{indep}.}  
Consider the quadratic map $Q:\R^n \to \R$ given by the formula $Q(V)=|V|^2$. The graph of $Q$ is a paraboloid. For $V\in \R^n$, let $\widetilde V \in \R^{n+1}$ be the vector $(V,Q(V))$. If $S$ is a simplex in $\R^n$, let $\widetilde S$ be its lift in $\R^{n+1}$.

Consider the linear map $\pi_i: \R^{n+1} \to \R^n$ given by the formula
$$
(x_1,\dots,x_{n+1}) \mapsto (x_1,\dots,x_{i-1},x_{n+1},x_{i+1},\dots,x_n).
$$
In terms of these maps, the $i$th coordinate of $C(S)$ is $\V(\pi_i(\widetilde S))/2\V(S)$. 

Let $F$ be a facet of $P$. Denote by $S_F$ the cone over $F$ with the vertex $O$. The 
statement that we need to prove is  that, for every $i$, 
\begin{equation} \label{summ}
\sum_{F\subset \partial P} \V(\pi_i(\widetilde S_F))
\end{equation}
is independent of $O$, where the sum is taken over the facets of $P$. Without loss of generality, consider the case $i=1$.

 Let $V_1,\dots,V_n$ be the vertices of a facet $F$. The respective summand is 
\begin{equation} \label{dett}
\det \left|\begin{array}{cccc}
|O|^2&|V_1|^2&\dots&|V_n|^2\\
\overline O&\overline V_1&\dots&\overline V_n\\
1&1&\dots&1
\end{array}\right| 
\end{equation}
where $\overline O=(o_2,\dots,o_n)$ and likewise for $\overline V_1,\dots,\overline V_n$.

The coefficients of $|O|^2,o_2,o_3,\dots,o_n$ in (\ref{dett}) are, up to the sign, 
$$
\V(p_1 \pi_1(\widetilde F)), \V(p_2 \pi_1(\widetilde F)), \dots, \V(p_n \pi_1(\widetilde F)),
$$ 
where the projection $p_i: \R^n \to \R^{n-1}$ is given by forgetting $i$th coordinate. 

It remains to notice that the total signed volume of $p_i \pi_1(\widetilde{\partial P})$ is equal to zero since $\partial P$ is a cycle. Therefore the
coefficients of $|O|^2,o_2,o_3,\dots,o_n$  in (\ref{summ}) vanish, and this sum does not depend on $O$, as needed.
\proofend

Similarly to Proposition \ref{coord}, one has an explicit formula for $CCM(P)$ that follows from Lemmas \ref{vol} and \ref{cc} and the proof of Theorem \ref{indep}. 

Let $F=(V_1,\dots,V_n)$ be a face of a simplicial polyhedron $P\subset \R^n$. Let $A(F)$ be the matrix whose columns are the vectors $V_1,\dots,V_n$, and let $A_i(F)$ be obtained from the matrix $A(F)$ by replacing its $i$th row by the row $(|V_1|^2,\dots,|V_n|^2)$. Choosing $O$ as the origin yields the following formula for the $i$th component of the Circumcenter of Mass.

\begin{proposition} \label{multicoord}
One has:
$$
CCM(P)_i=\frac{1}{2 (n!) V(P)} \sum_{F\subset \partial P} \det A_i(F).
$$
\end{proposition}

This rational function of the coordinates of the vertices of the polyhedron $P$ has homogeneous degree one: it is a  ratio of polynomials of degrees $n+1$ and $n$.

As an application of this formula, let us prove an analog of Corollary \ref{inscr}. 

\begin{corollary} \label{multidinscr}
The Circumcenter of Mass of an inscribed polyhedron is its circumcenter.
\end{corollary}

\proof
Choose the origin $O$ at the circumcenter and, without loss of generality, assume that the circumscribing sphere is unit. We wish to prove that $CCM(P)=O$; to be concrete, consider the $n$th component of $CCM(P)$, as given in Proposition \ref{multicoord}. 

The matrices $A_n(F)$ have $n$th row consisting of $1$s, that is, these are the matrices from Lemma \ref{vol}. It follows that, up to a factor, $CCM(P)_n$ is the signed volume of the projection of $\partial P$ to $\R^{n-1}$ along the last coordinate direction. This volume vanishes, as in the proof of Theorem \ref{indep}, and we are done. 
\proofend

A multi-dimensional version of Theorem \ref{Archimedes} holds as well, and for the same reason: the contribution of the simplices whose bases belong to the cut cancel out. Thus, with the self-explanatory notation, one has the Archimedes Lemma.

\begin{theorem} \label{ArchimedesForPolytopes}
If a simplicial polyhedron $P$ is decomposed into simplicial polyhedra $Q$ and $R$ then 
$CCM(P)=CCM(Q\oplus R).$
\end{theorem}

It follows that one can find the Circumcenter of Mass using any triangulation of $P$.

Similarly to the 2-dimensional case, one can take an affine combination of the centroid and the circumcenter of a simplex to define the center $C_t(P)$; as $t$ varies, these centers lie on a line that we call the Euler line of the polyhedron. 

The Euler line of a simplex is well studied: this is the line through the centroid and the circumcenter, see \cite{BB,HW} and the references in the latter article. 

The Euler line of a simplex contains its {\it Monge point}. This is the intersection point of the hyperplanes through the centroids of the $n-2$-dimensional faces of a simplex in $\R^n$, perpendicular to the opposite 1-dimensional edges (in dimension two, the Monge point is the orthocenter). In our notation, the Monge point is 
$C_{(n+1)/(n-1)}$. Thus we obtain a definition of the Monge point of a simplicial polyhedron in $\R^n$.

\section{On the sphere} \label{sphere}

In this section, we consider a version of the theory on the sphere (and, to some extent, in the hyperbolic space). Let us start with $S^2$. 

A circle on the unit sphere in $\R^3$ is its section by an affine plane, that is, a plane not necessarily through the origin which is the center of the sphere. To an oriented circle we assign its center; the choice of the two antipodal centers is made using the right-hand rule. Given a spherical triangle $ABC$, its circumcircle is oriented by the cyclic order of the vertices, and this determines the circumcenter of the oriented triangle.

We need the notion of the center of mass of a collection of points on the sphere, see \cite{Ga}. Let $V_i \in S^2,\ i=1,\dots,n$, be points and $m_i \geq 0$ their masses. The center of mass of this system of masses is the point
$$
\frac{\sum m_i V_i}{|\sum m_i V_i|},
$$
taken with the mass $|\sum m_i V_i|$. This notion of the center of mass satisfies natural axioms and is uniquely characterized by them \cite{Ga}. 

The notion of the center of mass naturally extends to spherical lamina. If $U$ is a spherical domain then its center of mass is the vector
$$
\int_U v\ dA,
$$
normalized to be unit, where $v\in S^2$ is the position vector of a point and $dA$ is the standard area form on the sphere. 

In particular, let $ABC$ be a spherical triangle. Then the center of mass of the triangular lamina is the vector
\begin{equation} \label{Brock}
A\times B\ \frac{d(A,B)}{\sin d(A,B)} + B\times C\ \frac{d(B,C)}{\sin d(B,C)} + C\times A\ \frac{d(C,A)}{\sin d(C,A)}, 
\end{equation}
normalized to be unit. Here $A\times B$ is the cross-product in $\R^3$, and $d(A,B)$ is the spherical distance between  points, that is, the side length of the triangle. The mass of the triangular lamina is the norm of the vector (\ref{Brock}). Formula (\ref{Brock}) can be found in \cite{Br,HJN,Mi}.

With these preparations, we can define the circumcenter of mass of a spherical polygon and study its properties, similarly to the Euclidean case. 

Let $P=(V_1,\dots, V_n)$ be a spherical polygon and let $W\in S^2$ be a point, distinct from the vertices of $P$. Consider the triangles $W V_i V_{i+1}$ where the index $i$ is cyclic. Let $O_i$ be the circumcenter of $i$th triangle, taken with the mass $m_i$, equal to the area of the plane triangle $W V_i V_{i+1}$. Define $CCM(P)$ to be the center of mass of the collections of points $(O_i,m_i)$. 

\begin{theorem} \label{indepS}
The point $CCM(P)$ does not depend on the choice of $W$.
\end{theorem}

\proof Let $ABC$ be a spherical triangle. Its circumcenter is the intersection of the sphere with the ray, perpendicular to the plane of the triangle and oriented according to the orientation of the plane given by the cyclic order of the vertices. Hence the circumcenter is the vector $(A-B) \times (A-C)$, normalized to be unit, that is,
$$
\frac{A\times B + B\times C + C\times A}{|A\times B + B\times C + C\times A|}.
$$
Note that the denominator is twice the area of the plane triangle $ABC$.

Using this formula for the triangles $W V_i V_{i+1}$ and taking the weighted sum, yields 
\begin{equation} \label{CCMS}
CCM(P) = \frac{\sum V_i \times V_{i+1}}{|\sum V_i \times V_{i+1}|},
\end{equation}
since $\sum (W\times V_i + V_{i+1}\times W) =0$. Thus $CCM(P)$ is independent of $W$.
\proofend

The Archimedes Lemma holds, the same way as in the Euclidean case. We also have an analog of Theorem \ref{equilat}. 

\begin{theorem} \label{equilatS}
If $P$ is an equilateral spherical polygon then its circumcenter of mass coincides with the center of mass of the lamina bounded by $P$, that is, $CCM(P)=CM(P)$.
\end{theorem}

\proof Let $\ell$ be the side length of $P$. Choose a point $W$ and triangulate the polygon as above. Let $(Q_i,m_i)$ be the center of mass of the triangular lamina $WV_i V_{i+1}$. Then $CM(P)$ is the center of mass of the system of points $(Q_i,m_i)$.  By formula (\ref{Brock}), this is
$$
\sum V_i \times V_{i+1} \frac{\ell}{\sin \ell} + \sum W \times V_i \frac{d(W,V_i)}{\sin d(W,V_i)} + \sum V_{i+1} \times W \frac{d(W,V_{i+1})}{\sin d(W,V_{i+1})},
$$
as always, normalized to a unit vector. 
The last two sums cancel each other and we obtain $CCM(P)$, as given by formula (\ref{CCMS}).
\proofend

An analog of Theorem \ref{cont} holds as well. Let $\gamma(t)$ be an oriented simple closed spherical curve. The continuous limit of the sum (\ref{CCMS}) is the vector-valued integral
$$
\int_{\gamma} \gamma \times \gamma'\ dt
$$
(note that it does not depend on the parameterization). Denote this integral by $CCM(\gamma)$, and let $CM(\gamma)$ be the center of mass of the homogeneous lamina bounded by $\gamma$.

\begin{theorem} \label{contS}
One has: $CCM(\gamma)=CM(\gamma)$.
\end{theorem}

\proof Let $U$ be the domain bounded by $\gamma$.  We claim that 
$$
\int_U v\ dA = \frac{1}{2} \int_{\gamma} \gamma \times \gamma'\ dt,
$$
where, as before, $v$ is the position vector of a point.

Let $\xi$ be a test vector. Then what we need to establish is the equality
$$
\int_U v\cdot\xi\ dA = \frac{1}{2} \int_{\partial U} (\gamma \times \gamma')\cdot\xi\ dt.
$$
Since $v$ is the unit normal vector to the sphere, the left hand side is the flux of the constant vector field $\xi$ through $U$. The right hand side equals
$$
-\frac{1}{2} \int_{\partial U} (\gamma \times \xi)\ d\gamma.
$$
The result will follow from the Stokes formula once we show that 
$$
-\frac{1}{2} {\rm curl} (\gamma\times \xi) = \xi.
$$ 
The last equality is easily verified by a direct computation, writing 
$\gamma=(x,y,z),\ \xi=(a,b,c)$, 
where $x,y,z$ are the Cartesian coordinates in $\R^3$ and $a,b,c$ are constants.
\proofend

The case of higher-dimensional sphere $S^n\subset \R^{n+1}$ is similar. The cross-product is replaced by the skew-linear  operation
$$
(V_1,\dots,V_n) \mapsto V_1\times\dots\times V_n
$$
defined by the equality
$$
\det (V_1,\dots,V_n,\xi) = (V_1\times\dots\times V_n)\cdot \xi
$$
for every test vector $\xi$. 

Given a simplicial spherical polyhedron $P$,  one triangulates it by choosing a point $W$ and takes the weighted sum of the circumcenters of the respective simplices. This results in the vector
\begin{equation} \label{multid}
\sum V_1\times\dots\times V_n,
\end{equation}
normalized to be unit, where the sum is taken over the facets $(V_1,\dots,V_n)$ of the polyhedron $P$. 

Finally, we make a brief comment on the hyperbolic case. In this case, the Euclidean space is replaced with the Minkowski space, the pseudo-Euclidean space of signature $(1,n)$, and the sphere $S^n$ by the pseudo-sphere, the
hyperboloid model of the hyperbolic space $H^n$. One can still define cross-product, and one can use formula (\ref{multid}) for the circumcenter of mass. The center of mass is also  defined in the hyperbolic case, see \cite{Ga}.

However, the geometrical interpretation is not as straightforward as in the spherical case. To fix ideas, consider the case $n=2$, the hyperbolic plane. 

One has three possibilities for the vector (\ref{multid}): it may be space-like, time-like, or null. In the first case, the line spanned by this vector intersects the upper sheet of the hyperboloid $z^2-x^2-y^2=1$, and we obtain a point of the hyperbolic plane. In the last case, we obtain a point on the circle at infinity, and in the second case, a point on the hyperboloid of one sheet $z^2-x^2-y^2=-1$ (outside of our model of $H^2$). 

This issue is already present in the case of a triangle. 
The intersection of the upper sheet of the hyperboloid $z^2-x^2-y^2=1$ with an affine plane may be a closed curve, disjoint from the null cone $z^2=x^2+y^2$; it may be a closed curve tangent to the null cone; or it may be an open curve intersecting the null cone.
A triangle is circumscribed by a curve of constant curvature,  the intersection of the hyperboloid $z^2-x^2-y^2=1$ with the plane of the triangle.  If the curvature is greater than one, this curve is a circle; if the curvature equals one, it is a horocycle; and if the curvature is less than one, the curve is an equidistant curve. These three cases correspond to the three relative positions of the plane and the null cone.

All three types of curves are represented by circles in the Poincar\'e disc model. In the first case, the circle lies inside the Poincar\'e disc, in the second it is tangent to the boundary, and in the third it intersects the boundary. Only in the first case the center of the circumcenter of a triangle is a point of the hyperbolic plane. 

\begin{remark}
{\rm The results on the circumcenter of mass in the Euclidean case can be obtained in the limit $R\to\infty$ from the spherical case, where $R$ is the radius of the sphere.}
\end{remark}

\begin{remark}
{\rm A proof that the Bicycle (Darboux) Transformation \cite{TT} preserves the CCM of triangles and quadrilaterals in spherical and hyperbolic geometry is entirely analogous to the one in Euclidean space. This fact and conservation of area under the Transformation allows to fully describe the dynamics of these polygons under the Bicycle Transformation.}
\end{remark}

\bigskip
{\bf Acknowledgments}. It is a pleasure to acknowledge interesting discussions with V. Adler, A. Akopyan, I. Alevi, Yu. Baryshnikov, B. Gr\"unbaum, D. Hatch, I. Rivin, O. Radko, A. Sossinsky, A. Veselov. This project originated during the program Summer@ICERM 2012; we are grateful to ICERM for support and hospitality. 
S. T. was partially supported by the NSF grant DMS-1105442.


\begin{thebibliography}{99}

\bibitem{Ad1} V. Adler, {\it Cutting of polygons}.  Funct. Anal. Appl. {\bf 27} (1993), 141--143.

\bibitem{Ad2} V. Adler, {\it Integrable deformations of a polygon. } Phys. D {\bf 87} (1995),  52--57. 

\bibitem{AM} T. Apostol, M.  Mnatsakanian, {\it Finding centroids the easy way.} Math Horizons,  September,
2000, 7--12.



\bibitem{Br} J.  Brock, {\it Centroid and inertia tensor of a spherical triangle.} National Technical Information Service, Naval
Postgraduate School, Monterey, California, 1974.

\bibitem{BB} M. Buba-Brzozowa, {\it The Monge point and the $3(n+1)$ point sphere of an $n$-simplex.} J. Geom. Graph. {\bf 9} (2005), 31--36.

\bibitem{DD} E. Danneels, N.  Dergiades, {\it A theorem on orthology centers.}
Forum Geom. {\bf 4} (2004), 135--141.


\bibitem{Ga} G. Galperin, {\it A concept of the mass center of a system of material points in the constant curvature spaces.} Comm. Math. Phys. {\bf 154} (1993),  63--84.

\bibitem{HW} H. Havlicek, G. Weiss, {\it Altitudes of a tetrahedron and traceless quadratic forms.} Amer. Math. Monthly {\bf 110} (2003), 679--693.

\bibitem{HJN} S. Heilman, A. Jagannath,  A. Naor, {\it Solution of the propeller conjecture in $\R^3$}, preprint arXiv 1112.2993v1

\bibitem{Ho} T. Hoffmann, {\it Discrete Hashimoto surfaces and a doubly discrete smoke-ring flow}. Discrete differential geometry, 95--115, Oberwolfach Semin., 38, Birkh\"auser, Basel, 2008. 

\bibitem{La} C.-A. Laisant, {\it  Th\'eorie et applications des \'equipollences}. Gauthier-Villars, Paris 1887, pp. 150--151.


\bibitem{Ma} M.  Mammana, {\it The maltitude construction in a convex noncyclic quadrilateral}.  Forum Geom. {\bf 12} (2012), 243-245.

\bibitem{Mi} G.  Minchin, {\it A treatise on statics, containing some of the fundamental propositions in electrostatics.} London, Longmans, Green \& Co. 1877, p. 259.

\bibitem{M} A. Myakishev, {\it On two remarkable lines related to a quadrilateral.} 
Forum Geom. {\bf 6} (2006), 289--295.

\bibitem{PSW} U. Pinkall, B. Springborn, S. Weissmann, {\it A new doubly discrete analogue of smoke ring flow and the real time simulation of fluid flow}. J. Phys. A {\bf 40} (2007),  12563--12576.

\bibitem{RT} O. Radko, E. Tsukerman, {\it The perpendicular bisector construction, the isoptic point and the Simson line of a quadrilateral}. Forum Geom. {\bf 12} (2012), 161-189.


\bibitem{TT} S. Tabachnikov, E. Tsukerman, {\it On the discrete bicycle transformation}. Preprint arXiv:1211.2345.

\bibitem{Ts} E. Tsukerman, {\it The perpendicular bisector construction in $n$-dimensional Euclidean and non-Euclidean geometries}. J. Classical Geom., in print. 




\end{thebibliography}
\end{document}